\documentclass[11pt]{amsart}

\usepackage{amsmath,amsthm,amsfonts,amssymb}
\usepackage[all,cmtip]{xy}
\usepackage{verbatim}
 \usepackage{latexsym}
\newtheorem{thm}{Theorem}
\newtheorem{prop}[thm]{Proposition}

\newtheorem{lemma}[thm]{Lemma}
\newtheorem{definition}[thm]{Definition}
 
  \newtheorem{eg}[thm]{Example}

\theoremstyle{remark}

\newtheorem{case'}{Case}

\numberwithin{equation}{section}
\numberwithin{thm}{section}

\renewcommand{\leq}{\leqslant}
\def\({\left(}
\def\){\right)}
\def\[{\left[}
\def\]{\right]}
\def\={\quad = \quad}
\def\+{\quad + \quad}

\def\R{{\mathbb{R}}}
\def\Z{{\mathbb{Z}}}

\def\C{{\mathbb{C}}}
\def\T{{\mathbb{T}}}

\def \a {\alpha}
\def \b {\beta}

\def \l {\lambda}

\def \L {\Lambda}

\begin{document}

\title{A Density Condition for Interpolation on the Heisenberg Group}


\author{Bradley Currey and
Azita Mayeli}

\date{\today}
 
\maketitle

 \begin{abstract} Let $N$ be the Heisenberg group. We consider left-invariant multiplicity free subspaces of $L^2(N)$. We prove a necessary and sufficient density condition in order that such subsspaces possess the interpolation property with respect to a class of discrete subsets of $N$ that includes the integer lattice. We exhibit a concrete example of a subspace that  has interpolation for the integer lattice, and we also prove a necessary and sufficient condition for shift invariant subspaces to possess a singly-generated orthonormal basis of translates.

   \end{abstract}



\section{introduction}\label{intro}    

Let $\mathcal H$ be a Hilbert space of continuous functions on a topological space $X$ for which point evaluation $f \mapsto f(x)$ is continuous, let $\Gamma$ be a countable discrete subset of $X$, and let $p$ be the restriction mapping $ f \mapsto f|_\Gamma$ on $\mathcal H$. For the present work, the sampling problem is as follows: describe those pairs $(\mathcal H, \Gamma)$ for which  $p$  is a constant multiple of an isometry of $\mathcal H$ into $\ell^2(\Gamma)$. If $p$ is surjective then we say that $(\mathcal H, \Gamma) $ has the interpolation property. Sampling and interpolation has been studied by various authors in various settings; three related examples are \cite{Pesen98},  \cite{F}, and \cite{FG}. In \cite{Pesen98}  the  author studies sampling on stratified Lie groups,  
while some of the results in  \cite{F} provide a characterization of left invariant sampling subspaces of $L^2(G)$ where $G$ is any locally compact unimodular Type I topological group, in terms of the notion of admissibility. As a consequence of the fundamental work of \cite{F}, it has been suspected that for the Heisenberg group, left-invariant sampling spaces cannot have the interpolation property with respect to lattice subgroups. The work of \cite{FG} studies more general (non-tight) sampling and includes ideas from both of the preceding articles. 

Here we consider the interpolation property for a class of quasi-lattices $\Gamma_{\a,\b}, \a , \b > 0$ in the Heisenberg group $N$ that we also consider in \cite{CM08}.  Let $\mathcal H$ be a left invariant subspace of $L^2(N)$ that is multiplicity free: the group Fourier transform of functions in $\mathcal H$ have rank at most one. In an explicit version of the group Plancherel transform,  the dual $\hat N$ is a.e. identified with $\R\setminus \{0\}$, the Plancherel measure on $\hat N$ becomes $d\mu = |\l|d\l$, and there is a measurable subset $E$ of $\R\setminus \{0\} $ such that $\mathcal H$ is naturally identified with $L^2(E\times \R)$. With this identification, a left translation system $\{T_\gamma\psi : \gamma \in \Gamma_{\a,\b}\}$ where $\psi\in \mathcal H$  becomes a {\it field} over $E$ of Gabor systems in $L^2(\R)$. We use this identification to show that $(\mathcal H, \Gamma_{\a,\b})$ has the interpolation property exactly when the $\mu(E) = 1/\a\b$. 

The paper is organized as follows:
after introducing some preliminaries, in Section \ref{group} we collect relevant results from \cite{F, FG} concerning admissibility and sampling for left invariant subspaces of $L^2(G)$, where $G$ is unimodular and  type I. The point is that admissibility is necessary for sampling, and (Theorem \ref{samp-charac}) a left invariant subspace is sampling if and only if it is admissible and its convolution projection generates a tight frame. We also recall the general fact that such a subspace has interpolation property if the afore-mentioned Parseval frame is actually orthonormal.  
In Section \ref{Heisenberg} we specialize to the Heisenberg group $N$ and direct our attention to the interpolation property for  multiplicity free subspaces with respect to the discrete subsets $\Gamma_{\a,\b}$. In Theorem \ref{interpolMF} we characterize such spaces that have the interpolation property, and in Example \ref{mainEG} we give a concrete example of a muliplicity free subspace of $L^2(N)$ that has the interpolation property with respect to the integer lattice in $N$. Finally, with Theorem \ref{ONcharacter} we prove a necessary and sufficient condition for any shift invariant subspace to have the interpolation property.


\section{Sampling spaces for unimodular groups }\label{group}



Let $G$ be a locally compact unimodular, topological group that is type I 
 and choose a Haar measure $dx$ on $G$.  For each $x\in G$ let $T_x$ be the unitary left translation operator on $L^2(G)$. Let $\hat G$ be the unitary dual of $G$, the set of equivalence classes of continuous unitary irreducible representations of $G$, endowed with the hull-kernel topology and the Plancherel measure $\mu$. 
As is well-known, for each $\l \in \hat G$, there is a continuous unitary irreducible representation $\pi_\l$ belonging to $\l$, acting in a Hilbert space $\mathcal L_\l$  with the following properties.

\vspace{.1in}
\noindent
(1) For each $\phi \in L^1(G) \cap L^2(G)$, the weak operator integral
$$
\pi_\l(\phi) := \int_N \phi(x) \pi_\l(x) dx
$$
defines a trace-class 
operator on $\mathcal L_\l$.

\vspace{.1in}
\noindent
(2) The 
group Fourier transform
$$
\mathcal F : L^1(G) \cap L^2(G) \rightarrow \int_{\hat G}^\oplus \ \mathcal{HS}(\mathcal L_\l) d\mu(\l)
$$
defined by $\phi \mapsto \{\pi_\l(\phi)\} := \{\hat\phi(\l)\}_{\l\in \hat G}$ satisfies $\| \mathcal F(\phi)\| = \|\phi\|_2$ and has dense range. (Here $\mathcal{HS}(\mathcal L_\l)$ denotes the Hilbert space of Hilbert-Schmidt operators on $\mathcal L_\l$.)

\vspace{.1in}
\noindent
(3) For each $x\in G$, $\mathcal F(T_x\phi) = \pi_\l(x) \hat\phi(\l), \l \in \hat G$.

\vspace{.1in}
\noindent
A closed subspace $\mathcal H$  of $L^2(G)$ is said to be  \it left invariant
 \rm if $T_x(\mathcal H) \subset \mathcal H$ holds for all $x \in G$.  Let $\mathcal H$ be a left invariant subspace of $L^2(G)$, and let $P : L^2(G) \rightarrow \mathcal H$ be the orthogonal projection onto $\mathcal H$. Then there is a unique (up to $\mu$-a.e. equality) measurable field $\{\hat P_\l\}_{\l\in \hat G}$ of orthogonal projections where $\hat P_\l$ is defined on $\mathcal L_\l$, and so that 
 $$
 \widehat{(P\phi)}(\l) = \hat\phi(\l) \hat P_\l
 $$
 holds for $\mu$-a.e.  $\l \in \hat G$. Set $m_\mathcal H(\l) = \text{rank}(\hat P_\l)$. Then spectrum of $\mathcal H$ is the set $\Sigma(\mathcal H) = \text{supp}(m_\mathcal H)$.  
  A left invariant subspace $\mathcal H$ of $L^2(G)$ is said to be multiplicity free if  $m_\mathcal H(\l) \le 1$ a.e.; if $\mathcal H$ is left invariant and  $m_\mathcal H(\l) = 1$ a.e. then we will say that  $\mathcal H$ is multiplicity one.  
   
    Recall that $\psi \in \mathcal H$ is said to be admissible (with respect to the left regular representation) if the operator $V_\psi$ defined by $V_\psi(\phi) = \phi * \psi^*$ defines an isometry of $\mathcal H$ into $L^2(G)$. For convenience we recall   \cite[Theorem 4.22]{F}  for unimodular groups.    
    
  \begin{thm}\label{F-4.22}
 Let $\mathcal H$ be a  closed left invariant subspace of $L^2(G)$
with the associated  measurable projection field $\{\hat P_\l \}$. Then $\mathcal H$ has   admissible vectors if and  only if the map $\l\mapsto \text{rank}(\hat P_\l)$ is finite and  $\mu$-integrable. 
\end{thm} 
For an example of $\mathcal H$ and  an admissible vector we refer the interested reader to \cite{M08}. Different examples  have also been presented in this work.     
    Gleaning from results in \cite{F}, we have the following.
 
 \begin{prop} \label{admissible} Let $\mathcal H$ be a closed left invariant subspace of $L^2(G)$ with $G$ unimodular. Then the following are equivalent.
 
 \vspace{.1in}
 \noindent
 (a)  $\mathcal H$ has an admissible vector. 
 
 \vspace{.1in}
 \noindent
 (b) There is a left invariant subspace $\mathcal K$ of $L^2(G)$ and  $\eta \in\mathcal K$ such that $\phi \mapsto \phi * \eta^*$ is an isometric isomorphism of $\mathcal K$ onto $\mathcal H$. 
 
   \vspace{.1in}
 \noindent
 (c) There is a unique self-adjoint  convolution idempotent $S\in \mathcal H$ such that $\mathcal H = L^2(G) * S$. 
 
  \vspace{.1in}
 \noindent
 (d) The function $m_\mathcal H$ is integrable over $\hat G$ with respect to Plancherel measure $\mu$.

 \end{prop}
 
  \begin{proof} Let  (a) hold and $\psi$ be an admissible vector for $\mathcal H$. Then $V_\psi^* V_\psi$ is an isometry and hence the identity on $\mathcal H$.  Take 
    $\mathcal K = V_\psi(\mathcal H)$. Then $V_\psi^*$ is an isometric isomorphism of  $\mathcal K$ onto $\mathcal H$. To prove  (b) , we only need to show that $V_\psi^*$ acts by $V_\psi^* f= f\ast \psi$. For this, let $f\in \mathcal H\ast \psi$ and $g\in L^1(G)\cap L^2(G)$. Then 
    \begin{align}\label{convolution-op}
    \langle V_\psi^\ast f, g\rangle = \langle f, V_\psi g\rangle= \langle f, g\ast \psi^\ast \rangle= \langle f\ast \psi, g\rangle.
    \end{align}
    Now   $\eta = \psi^*$ applies for (b).

   Suppose that (b) holds. Define $V_\eta(\phi)=\phi\ast \eta^*$ for $\phi\in \mathcal K$. Then  $V_\eta V_\eta^*$ is projection of $L^2(N)$ onto $\mathcal H$.
   Since $V_\eta^*$ is bounded,  then by an analogous computation in (\ref{convolution-op}) we have $V_\eta^* = V_{\eta^*}$, and  hence    $V_\eta V_\eta^*(\phi) = \phi  * (\eta * \eta^*)$.  Now set  $S = \eta*\eta^* = V_\eta(\eta)$. Evidently, $S$ belongs to $ \mathcal H$,  is self-adjoint, and is a convolution idempotent, and the projection  onto $\mathcal H$ is given by convolution with $S$.
 
 Suppose that (c) holds. Then $S$ itself is an admissible vector in $\mathcal H$, so (d) follows from Theorem \ref{F-4.22}. Finally,  Theorem \ref{F-4.22} says that (d) implies (a).
 
 \end{proof}

We will say that a left invariant subspace $\mathcal H$  is {\it admissible} if it satisfies one of  the conditions of Proposition \ref{admissible}. The function $S$ of condition (c)   is called the reproducing kernel for $\mathcal H$ . Note that in this case the associated projection field $\{\hat P_\l\}$ is just the group Fourier transform of $S$. 
      
 \begin{definition}
 Let $\Gamma$ be any countable discrete subset of $G$ and let $\mathcal{H}$  be a left invariant subspace of $L^2(G)$ consisting of continuous functions. We shall call  $(\mathcal{H},\Gamma)$ a sampling pair if 
  there exist $S\in \mathcal{H}$  and $c = c_{\mathcal H,\Gamma} > 0$ such that  for all $\phi\in  \mathcal{H}$
  

 
  \begin{align}\label{isometry}
   \| \phi\|^2= \frac{1}{c}  \sum_{\gamma\in \Gamma} | \phi(\gamma)|^2,
 \end{align}
   and
 \begin{align}\label{sinc-equality}
 \phi= \sum_\gamma \phi(\gamma)~ T_\gamma S.
 \end{align}
 where the sum (\ref{sinc-equality}) converges in $L^2$. 
  \end{definition}
  
Following \cite{F}  we say that $S$ is a sinc-type function. It is well-known that the sum (\ref{sinc-equality}) above converges uniformly as well as in $L^2$ (see \cite[Remark 2.4]{FG}).

Recall that a system $\{\psi_j\}_{j\in J}$ of functions in a separable Hilbert space $\mathcal H$ is a tight frame for $\mathcal H$ if for some $c > 0$, 
$$
c \ \|g\|^2 = \sum_{j\in J} \ \left|\langle g, \psi_j\rangle \right|^2
$$
holds for every $g \in \mathcal H$. A Parseval frame is a tight frame for which $c = 1$. Now suppose that  $\mathcal H$ is closed left invariant admissible with reproducing kernel $S$, and that $\Gamma$ is a countable discrete subset such that the relation  (\ref{isometry}) holds for all $\phi\in\mathcal H$. 
Then the identity $\phi(x)  = \phi * S(x) = \langle \phi, T_x \rangle$ makes it clear that $\T_\gamma S : \gamma \in \Gamma\}$ is a tight frame with frame bound $c$ and hence that $\mathcal H $ is a sampling space with sinc-type function $\frac{1}{c} S$ (see [Corollary 2.3] \cite{FG}.

To characterize sampling spaces, we only need to  observe that every such subspace is necessarily admissible. 


\begin{thm}\label{bdd-spect}
Let $\mathcal{H}$ be a left invariant subspace consisting of continuous functions and suppose that for some countable discrete subset $\Gamma$ of $G$,  (\ref{isometry})  holds for all $\phi \in \mathcal H$. Then  $\mathcal{H}$ is   admissible and hence a  sampling space. 
\end{thm}

\begin{proof} This is an immediate consequence of \cite[Theorem 2.56]{F}; see also \cite[Theorem 2.2]{FG}. 
\end{proof}

Hence we have the following equivalent conditions for left invariant subspaces.

\begin{thm}\label{samp-charac}
Let $\mathcal{H}$ be a left invariant subspace and let $\Gamma$ be a countable discrete subset of $G$. Then the following are equivalent.

\vspace{.1in}
\noindent
(i)  $\mathcal H$ is admissible and $T_\gamma S : \gamma\in \Gamma\}$ is a tight frame for $\mathcal H$ with frame bound $c$, where $S$ is its reproducing kernel. 

\vspace{.1in}
\noindent
(ii)  $\mathcal H$ consists of continuous functions and the map $A: \mathcal{H}\rightarrow \ell^2(\Gamma)$ defined by $A(\phi)= \{\frac{1}{\sqrt{c}}\phi(\gamma)\}$ is isometry.
\end{thm}

Of course if $\mathcal H$ satisfies one of the conditions of the preceding theorem, then  $(\mathcal{H},\Gamma)$ is a sampling pair with sinc-type function $\frac{1}{c}S$. Next we turn to the question of interpolation. 

\begin{definition} 
We say a  sampling pair $(\mathcal{H}, \Gamma)$ has an interpolation property if the isometry map $A$ is also  surjective. 
\end{definition}

The following is a consequence of Theorem \ref{samp-charac} and standard frame theory.

\begin{thm}\label{surj-onb}
Let    $(\mathcal{H}, \Gamma)$ be a sampling pair. Then there exists  a  sinc-type function $S$ for $\mathcal H$ for which  the following equivalent  properties hold:   $(\mathcal{H}, \Gamma)$  has the interpolation property if and only if  $\{\frac{1}{\sqrt{c }}T_\gamma S : \gamma \in \Gamma\}$ is an orthonormal basis for $\mathcal{H}$. 
\end{thm}

\begin{proof} By Theorem \ref{samp-charac}, $\mathcal H$ is admissible, and denoting its convolution projection by $S$,  we have that $\{\frac{1}{\sqrt{c}}T_\gamma S\}_{\gamma \in \Gamma}$ is a Parseval frame for $\mathcal H$ and the isometry $A$ is the associated analysis operator. If $\{\frac{1}{\sqrt{c}}T_\gamma S\}_{\gamma \in \Gamma}$ is an orthonormal basis then of course $A$ is surjective. On the other hand, if the isometry $A$ is surjective then $A$ is unitary and if $\delta_\gamma$ denotes the canonical basis element in $\ell^2(\Gamma)$ then $\| \frac{1}{\sqrt{c}}T_\gamma S \| = \| A^*\delta_\gamma\| = \| \delta_\gamma\| = 1$.

\end{proof}

In the following section we describe a class of subspaces of the Heisenberg group that admit sampling with the interpolation property. 


\section{The Heisenberg group and multiplicity free subspaces }\label{Heisenberg}

We now assume that $G=N$ is the Heisenberg group: as a topological space $N$ is identified with $\R^3$, and we let $N$ have the group operation
$$
(x_1,x_2,x_3)\cdot (y_1,y_2,y_3) = (x_1 + y_2,x_2+y_2,x_3 + y_3+x_1y_2).
$$
We recall some basic facts about harmonic analysis on $N$. Put $\L = \R \setminus \{0\}$. 
For $x \in N$, $\l \in \L$,  we define the unitary operator $\pi_\l(x)$ on $L^2(\R)$ by
$$
\Bigl(\pi_\l(x)f\Bigr)(t) = e^{2\pi i \l x_3} e^{-2\pi i \l x_2 t} f(t-x_1), \ f \in   L^2(\R).
$$
Then $x \mapsto \pi_\l(x)$ is an irreducible representation of $N$ (the Schr\"odinger representation), and for $\l \ne \l'$, the representations $\pi_\l$ and $\pi_{\l'}$ are inequivalent. With respect to the Plancherel measure on $\hat N$, almost every member of $\hat N$ is realized as above and the group Fourier transform takes the explicit form
$$
\mathcal F : L^2(N) \rightarrow \int_\L^\oplus \mathcal{HS}\bigl(L^2(\R)\bigr) |\l| d\l.
$$
Let $\mathcal H$ be a multiplicity free subspace of $L^2(N)$, $E = \Sigma(\mathcal H)$ the spectrum of $\mathcal H$, $P$ the projection onto $\mathcal H$, and let $\{\hat P_\l\}$ be the associated measurable field of projections. We have a measurable field $e = \{e_\l\}_{\l\in \L}$ where each  $e_\l$ belongs to  $L^2(\R)$, where $(\l \mapsto \|e_\l\|) = \bold 1_E$, and where $\hat P_\l = e_\l \otimes e_\l$ (i.e.  $\mathcal K_\l = \C e_\l$) for $\l\in E$. Thus the image of $\mathcal H$ under the group Fourier transform is
\begin{equation}\label{mult one int}
\hat{\mathcal H}= \int^\oplus_E \ L^2(\R) \otimes e_\l \ |\l| d\l.
\end{equation}
where $e_\l$ is regarded as an element of $\overline{L^2(\R)}$.
Hence $\mathcal H$ is isomorphic with
\begin{align}\label{reduced-subspace}
\int^\oplus_E \ L^2(\R) \ |\l|d\l
\end{align}
via the unitary isomorphism $V_e$ defined on $\mathcal H$ by $\{V_e\eta(\l) \}_{\l \in E}$ where $\eta \in \mathcal H$ and 
$$
V_e\eta(\l) = \hat\eta(\l)\bigl(e_\l\bigr), \ {\rm a.e.\ } ~ \l\in E.
$$
We identify the direct integral (\ref{reduced-subspace})  with $L^2(E\times\R)$ in the obvious way, where it is understood that $E$ carries the measure $|\l|d\l$. 
Note that if we write $\hat\eta(\l) = \{ f_\l\otimes e_\l\}_\l$, then  $V_e\eta(\l) = f_\l$. For a fixed unit vector field $e=\{e_\l\}$ we will say that $V_e$ is the reducing isomorphism for $\mathcal H$ associated with the vector field $e$. Note that given a multiplicity-free  subspace $\mathcal H$, the unit vector field $e =\{e_\l\}$ is essentially unique:  if $e' = \{e'_\l\}$ is another measurable unit vector field for which (\ref{mult one int}) holds, then there is a measurable unitary complex-valued function $c(\l)$ on $E$ such that $e'_\l = c(\l)e_\l$ holds for a.e. $\l$. Finally, given any subset $E$ of $\L$ and  measurable field $e = \{e_\l\}_{\l\in \L}$ with $(\l\mapsto \|e_\l\|) = \bold 1_E$, the subspace 
$$
\mathcal{H}_e = \{ \phi \in L^2(N):~ \text {Range}(\hat \phi(\l)^\ast)\subset \C e_\l , \text{ a.e. } \l\}
$$
is multiplicity free with spectrum $E$ and associated vector field $e$.

Let $\Gamma$ be a countable discrete subset of $N$.  If $V_e : \mathcal H \rightarrow L^2(E\times \R)$ is a reducing isomorphism, and $\psi \in \mathcal H$ with $g = V_e\psi$, then the system $\mathcal T(\psi, \Gamma) = \{T_{(k,l,m)} \psi : (k,l,m)\in\Gamma\}$ is obviously equivalent with the system  $ \widehat{ \mathcal  T}(g, \Gamma) = \{\hat T_{k,l,m}g : (k,l,m)\in\Gamma\}$ through the isomorphism $V_e$,  where 
$$
\hat T_{k,l,m}g(\l,t) = e^{2\pi i \l m}e^{-2 \pi i \l l t} \ g(\l,t- k).
$$
In \cite{CM08}, the discrete subsets $\Gamma_{\a,\b} = \a \Bbb Z \times \b\Bbb Z \times \Bbb Z$ for positive integers $\a$ and $\b$ are considered and when $\Gamma = \Gamma_{\a,\b}$, then we denote the above function systems by $\mathcal T(\psi, \a,\b)$ and  $ \widehat{ \mathcal  T}(g, \a,\b)$, respectively. For $\l \in \L$ fixed, $\hat T_{k,l,0}$ defines a unitary (Gabor) operator
on $L^2(\R)$ in the obvious way which we denote by $\hat
T_{k,l}^\l$. For $u \in L^2(\R)$ set   $\mathcal G(u,\a, \b,\l) = \{\hat
T_{k,l}^\l u : (k,l,0) \in \Gamma_{\a,\b}\}$. We say that $g \in L^2(E\times \R)$ is a Gabor field over $E$  with respect to $\Gamma_{\a,\b}$  if, for a.e. $\l \in E$, $\mathcal G(|\l|^{1/2} g(\l,\cdot),\a,\b,\l) $ is a Parseval frame  for $L^2(\R)$. If $g$ is a Gabor field over $E$ with respect to $\Gamma_{\a,\b}$, then standard Gabor theory implies that $\| |\l|^{1/2} g(\l,\cdot) \|^2 = \a\b|\l| \le 1$.

The following is an easy but significant extension of part of  \cite[Proposition 2.3]{CM08}. 



\begin{prop} \label{gabor fld} Let $E$ be a measurable subset of $\L$ and $g \in L^2(E\times \R)$ such that $ \widehat{ \mathcal  T}(g, \a,\b)$ is a Parseval frame for $L^2(E\times \R)$. Then $g$ is a Gabor field over $E$ with respect to $\Gamma_{\a,\b}$. 

\end{prop}

\begin{proof} Write $E = \cup E_j$ where each $E_j$ is translation congruent with a subset of $[0,1]$ and let $g_j = g|_{E_j}$. Then for each $j$, $ \widehat{ \mathcal  T}(g_j, \a,\b)$ is a Parseval frame for $L^2(E_j\times \R)$. Hence by \cite[Proposition 2.3]{CM08}, $g_j$ is a Gabor field over $E_j$, hence $g$ is a Gabor field over $E$. 

\end{proof}

The ``only part" of \cite[Proposition 2.3]{CM08} is also true for any subset $E$ provided the orthogonality of some  coefficient operators holds.

\begin{prop}
Let  $E\subset \L$ and   $g\in L^2(E\times \R)$ such that 
$g$ is 
   a Gabor field for $L^2(\R)$ over $E$ for $\Gamma_{\a,\b}$.  Write  $E = \dot\cup_{j\in J} E_j$ where each $E_j$ is translation congruent with a subset of $[0,1]$ and let $g_j = g{\bf |}_{E_j}$.
 For any $j$, let $C_j$ be the coefficient operator defined from $L^2(E\times \R)$  into $\ell^2(\Gamma_{\a,\b})$  by 
\begin{align}
C_j:~ f\rightarrow \{\langle \hat T_\gamma g_j, f\rangle \}_\gamma
\end{align}
and assume that for each $j \ne j'$, Range$(C_j) \subset $ Range$(C_{j'})^\perp$. Then $ \widehat{ \mathcal  T}(g, \a,\b)$ is a Parseval frame for $L^2(E\times \R)$. 
\end{prop}

\begin{proof} 
By \cite[Proposition 2.3]{CM08},
for any $j\in \Z$ the system $ \widehat{ \mathcal  T}(g_j, \a,\b)$ is a Parseval frame for $L^2(E_j\times \R)$. Therefore, with the orthogonality assumption, for any $f\in L^2(E\times \R)$  one has 
\begin{align}
\sum_\gamma \mid \langle \hat T_\gamma g, f\rangle\mid^2 =  \sum_\gamma \mid \sum_j \langle \hat T_\gamma g_j, f \rangle\mid^2 
 =\sum_{j} \sum_\gamma \mid \langle \hat T_\gamma g_j ,  f  \rangle  \mid^2   
= \parallel f\parallel^2,
\end{align}
and hence 
$ \widehat{ \mathcal  T}(g, \a,\b)$ is a Parseval frame for $L^2(E\times \R)$.
\end{proof}



The above observations together with Theorem \ref{samp-charac} above now give the following.

   \begin{thm}\label{e-sinc}  Let $E$ be a measurable subset of $\L$  and let $e = \{e_\l\}_{\l\in \L}$ be a measurable field of unit  vectors in $L^2(\R)$ such that $(\l\mapsto \|e_\l\|) = \bold 1_E$. The following are equivalent.
   
   \vspace{.1in}
\noindent
(i)  $E$ has finite Plancherel measure and $\widehat{ \mathcal  T}(\frac{1}{\sqrt{c}}e, \a,\b)$ is a Parseval frame for $L^2(E\times \R)$.

\vspace{.1in}
\noindent
(ii)  $(\mathcal H_e, \Gamma_{\a,\b})$ is a sampling pair with the sinc-type function $S=\frac{1}{c}V_e^{-1}(e)$.




\vspace{.1in}
\noindent
Moreover, if the above conditions hold, then $\frac{1}{\sqrt{c}}e$ is a Gabor field over $E$, and $E$ is included in the interval $[-1/\a\b,1/\a\b]$. 


   \end{thm} 
   
We now have a precise density criterion for the interpolation property in this situation.

\begin{thm} \label{interpolMF} Let $\mathcal H$ be a multiplicity free subspace of $L^2(N)$ with $E = \Sigma(\mathcal H)$. Suppose that for some $\a, \b > 0$, $(\mathcal H , \Gamma_{\a,\b})$ is a sampling pair with $c = c_{\mathcal H, \Gamma_{\a,\b}}$. Then $c =1/ \a\b$. Moreover,  $(\mathcal H , \Gamma_{\a,\b})$  has the interpolation property if and only if  $\mu(E)= 1 / \a\b$. Hence if  $(\mathcal H , \Gamma_{\a,\b})$  has the interpolation property, then $\a\b \le 1$.

\end{thm}

\begin{proof} By Theorem \ref{samp-charac}, $\mathcal H$ is admissible;  let $S$ be the associated reproducing kernel, and let $V_e$ be a reducing isomorphism, where $e = \{e_\l\}$ is an $L^2(\R)$-vector field with $(\l \mapsto \|e_\l \|) = \bold 1_E$, so that $V_e(S)=e$. It follows from the above that  $\{\frac{1}{\sqrt{c}} T_\gamma S \}_\gamma$ is a Parseval frame for $\mathcal H$, $\widehat{ \mathcal  T}(\frac{1}{\sqrt{c}}e, \a,\b)$ is a Parseval frame for $L^2(E\times \R)$, and  $\frac{1}{\sqrt{c}} e$ is a Gabor field over $E$. Hence for a.e. $\l \in E$, we have
$$
\||\l|^{1/2} \frac{1}{\sqrt{c}} e_\l \|^2 = \a\b |\l|, 
$$
and the relation $c = 1/ \a\b$  follows immediately. Now $\mathcal H$ has the interpolation property if and only if $\{\frac{1}{\sqrt{c}}T_\gamma S : \gamma \in \Gamma_{\a,\b}\}$ is an orthonormal basis for $\mathcal H$, if and only if $ \|\frac{1}{\sqrt{c}}S\|^2  = 1$. But 
$$
 \|\frac{1}{\sqrt{c}}S\|^2  =\a\b \|S\|^2 = \a\b \|V_e(S)\|^2 = \a\b \int_E \ \|e_\l\|^2 |\l| d\l =  \a\b \int_E |\l| d\l=\a\b\mu(E).
$$

This proves the first part of the theorem. Now if  $(\mathcal H , \Gamma_{\a,\b})$  has the interpolation property, then, since $E \subseteq [-1/\a\b,1/\a\b]$, we have  
$$ 1/\a\b = \int_E |\l| d\l \le \int_{[-1/\a\b,1/\a\b]} |\l| d\l = 1/(\a\b)^2.
$$



\end{proof}
 
We now construct an example of a sampling pair with the interpolation property. We assume that $\a = \b = 1$; note that in this case the interpolation property is equivalent with $E = [-1,1]$. In light of Theorems \ref{e-sinc} and \ref{interpolMF}, it is evident that in order to construct an example $(\mathcal H, \Gamma_{1,1})$ with $\mathcal H$ multiplicity free, it is enough to construct a measurable field of $L^2(\R)$-vectors  $\{e_\l\}$   such that $(\l \rightarrow \|e_\l \|)={\bf 1}_{[-1,1]}$ and such that $e$ generates a Heisenberg frame for $L^2([-1,1]\times \R)$. The  following technical lemma is helpful in the construction of  such a function $e$.




\begin{lemma} \label{tech}
Let  $e \in L^2([-1,1]\times \R)$ such that $e$ is a Gabor field over $[-1,1]$ with respect  to $\Gamma_{\a,\b}$, and such that the orthogonality condition 
\begin{equation}\label{orthog}
\sum_{k,l}\  \langle f(\l-1,\cdot), e_{k,l,0}(\l-1,\cdot)\rangle \ \overline{\langle f(\l,\cdot),e_{k,l,0}(\l,\cdot)\rangle} = 0
\end{equation}
 holds for all $\l\in (0,1]$ and for all  $f \in L^2([-1,1]\times \R)$.  Then the system $\widehat {\mathcal T}(e, \Gamma_{\a,\b})$ is a Parseval frame for $L^2([-1,1]\times \R)$. 

\end{lemma}

\begin{proof}
Suppose that $e$ is a Gabor field satisfying (\ref{orthog}) and let $f \in L^2([-1,1]\times \R)$. By Proposition 2.3 of \cite{CM08}, and the Parseval identity for Fourier series, we have 
$$
\begin{aligned}
\int_0^1 \|f(\l-1,\cdot)\|^2 |\l-1| d\l &= \sum_{k,l,m} \left| \int_0^1 \langle f(\l-1,\cdot),e_{k,l,m}(\l-1,\cdot)\rangle |\l-1| d\l \right|^2 \\
&= \sum_{k,l} \int_0^1 \left| \langle f(\l-1,\cdot),e_{k,l,0}(\l-1,\cdot)\rangle |\l-1| \right|^2 d\l 
\end{aligned}
$$
and similarly,
$$
\int_0^1 \|f(\l,\cdot)\|^2 |\l| d\l = \sum_{k,l} \int_0^1 \left| \langle f(\l,\cdot),e_{k,l,0}(\l,\cdot)\rangle |\l| \right|^2 d\l. 
$$
Hence
\begin{equation} \label{fnorm}
\begin{aligned}
\|f\|^2 &= \int_0^1 \|f(\l-1,\cdot)\|^2 |\l-1| d\l + \int_0^1 \|f(\l,\cdot)\|^2 |\l| d\l \\
&= \int_0^1 \ \sum_{k,l} \left( \bigl| \langle f(\l-1,\cdot),e_{k,l,0}(\l-1,\cdot)\rangle |\l-1| \bigr|^2 +  \bigl| \langle f(\l,\cdot),e_{k,l,0}(\l,\cdot)\rangle |\l| \bigr|^2 \right) \ d\l.
\end{aligned}
\end{equation}
But (\ref{orthog}) implies that for $\l \in (0,1]$, 
$$
\begin{aligned}
\sum_{k,l} & \left( \bigl| \langle f(\l-1,\cdot),e_{k,l,0}(\l-1,\cdot)\rangle |\l-1| \bigr|^2 +  \bigl| \langle f(\l,\cdot),e_{k,l,0}(\l,\cdot)\rangle |\l| \bigr|^2\right) \\
 &= 
\sum_{k,l} \left(\Bigl| \langle f(\l-1,\cdot),e_{k,l,0}(\l-1,\cdot)\rangle |\l-1| +   \langle f(\l,\cdot),e_{k,l,0}(\l,\cdot)\rangle |\l| \Bigr|^2\right).
\end{aligned}
$$
Combining the preceding with (\ref{fnorm}) and applying the Parseval identity for Fourier series again, we have
$$
\begin{aligned}
\|f\|^2 &= \sum_{k,l} \ \int_0^1 \ \Bigl| \langle f(\l-1,\cdot),e_{k,l,0}(\l-1,\cdot)\rangle |\l-1| + \langle f(\l,\cdot),e_{k,l,0}(\l,\cdot)\rangle |\l| \Bigr|^2 d\l \\
&= \sum_{k,l} \sum_m \ \Bigl|  \int_0^1 \left( \langle f(\l-1,\cdot),e_{k,l,0}(\l-1,\cdot)\rangle |\l-1| +  \langle f(\l,\cdot),e_{k,l,0}(\l,\cdot)\rangle |\l| \right) e^{-2\pi i\l m } d\l \Bigr|^2 \\
&= \sum_{k,l,m} \ \left| \int_{-1}^1 \langle f(\l,\cdot),e_{k,l,m}(\l,\cdot)\rangle |\l| d\l\right|^2.
\end{aligned}
$$
This proves the claim. 
\end{proof}

By virtue of Lemma \ref{tech}, it is sufficient  to construct a function $e$ with the properties in the preceding  lemma.
\begin{eg}\label{mainEG}
 For $\l \in (0,1]$, put 
$$
e_\l = \bold 1_{\left[\frac{1}{\l} - 1, \frac{1}{\l}\right]} \ \ \text{ and } \ \ \ e_{\l-1} = \bold 1_{[-1,0]}.
$$
Then $e$ defined by $e(\l, t)=e_\l(t)$ for  $\l\in (0,1]$ and $e(\l,t)= \bold 1_{[-1,0]}(t)$ for $\l\in [-1,0)$     is a Gabor field over $[-1,1]$ with respect to $\Gamma_{1,1}$. 

\end{eg}

\begin{proof} We compute that for any $f \in L^2([-1,1]\times\R)$ and for $\l \in (0,1]$, 
$$
\begin{aligned}
\langle f(\l-1,\cdot),e_{k,l,0}(\l-1,\cdot)\rangle  &= \int_\R f(\l-1, t) e^{2\pi i (\l-1) l t}\bold 1_{[-1,0]}(t-k) dt \\ &= \int_{I_k^{\l-1}} \ \left( \left(\frac{1}{1-\l}\right) \  f\left(\l-1, \frac{s}{\l-1}\right) \right) e^{2\pi i l s} ds
\end{aligned}
$$
and similarly, 
$$
\begin{aligned}
\langle f(\l,\cdot),e_{k,l,0}(\l,\cdot)\rangle  &= \int_\R f(\l,\cdot) e^{2\pi i \l l t}\bold 1_{\left[\frac{1}{\l} - 1, \frac{1}{\l}\right]} (t-k) dt \\ &= \int_{I_k^{\l}} \ \left( \frac{1}{\l} f\left((\l, \frac{s}{\l}\right)) \right) e^{2\pi i l s} ds
\end{aligned}
$$
where $I_k^{\l-1} = [-(1-\l)k, -(1-\l) k + (1-\l)]$ and $I_k^\l = [1 + \l k - \l, 1+ \l k]$. It is easily seen that for each $k$, 
$$
I_k^{\l-1} \cap I_k^\l = \emptyset \ \ \text{ and } \ \ (I_k^{\l-1} + k) \cup I_k^\l = [\l k , \l k + 1].
$$
Hence for each $k$,  the sequences $\{  \langle f(\l-1,\cdot), e_{k,l,0}(\l-1,\cdot)\rangle : l \in \Bbb Z\}$ and $\{\langle f(\l,\cdot),e_{k,l,0}(\l,\cdot)\rangle : l \in \Bbb Z \}$ are Fourier coefficients for orthogonal functions and we have
$$
\sum_l\  \langle f(\l-1,\cdot), e_{k,l,0}(\l-1,\cdot)\rangle \ \overline{\langle f(\l,\cdot),e_{k,l,0}(\l,\cdot)\rangle} = 0.
$$
Thus the equation (\ref{orthog}) holds for $e$. 
 
\end{proof}

Since the vector field $e = \{e_\l\}$ is compactly supported, one does not expect that the inverse Fourier image is well localized. We show this explicitly in the following, where we compute the inverse group Fourier transform in terms of ordinary Fourier transforms. For a function $f \in L^1(\R)$, put $\hat f(s) = \int_\R f(t) e^{2\pi i st} dt$ and $\check f(s) = \int_\R f(t) e^{-2\pi i st} dt$.

\begin{eg}    Let $ e = \{e_\l\}$ be the unit vector field from the preceding example, and let  $S \in L^2(N)$ be the function for which $V_e(S) = e$. For each $x\in\R$ define the intervals $I_{x,\l}$ and $J_x$ by 
$$
I_{x,\l} = \left[-\frac{1}{\l} - 1 ,\frac{1}{\l}\right] \bigcap \left(\left[-\frac{1}{\l} - 1 ,\frac{1}{\l}\right] + x\right),  \ J_{x} = [-1,0] \cap ([-1,0] + x).
$$
Then $S = S_0 + S_1$, where $S_0(x,y,z) = \check{F}_{x,y}(z)$ and $S_1(x,y,z) = \check{G}_{x,y}(z) $, and
where $G_{x,y}(\l) =\l   \bold 1_{[0,1]}(\l) \ \widehat{\bold 1}_{I_{x,\l}}(\l y)$ and $F_{x,y}(\l) = -\l  \bold 1_{[-1,0]}(\l) \ \widehat{\bold 1}_{J_{x}}(\l y)$. In particular, $S$  vanishes outside the strip $U = \{(x,y,z) : |x| < 1\}$ and $S_0$ and $S_1$ are given by sinc-type expressions. For example, for $(x,y,z)\in U$ and $y \ne 0, z \ne 0$, then 
$$
S_0(x,y,z) = \begin{cases} 
 \frac{1}{2\pi i   y} \left( \frac{e^{2\pi i (z-xy)}-1}{2\pi i (z-xy)} -
   \frac{e^{2\pi i (z+y)}-1}{2\pi i (z+y)}\right),
   & \text{ if} \ -1<x<0, ~z\neq xy, ~y\neq -z,\\
    \frac{1}{2\pi i   y}  \left(  \frac{e^{2\pi i z}-1}{2\pi i z}  -  \frac{e^{2\pi i (z+y(1-x))}-1}{2\pi i (z+y(1-x))}  \right), & \text{ if}\  \ 0<x< 1,~ x\neq -y(1-x).
    \end{cases}
$$

\end{eg}

\begin{proof}
 We have
$$
\begin{aligned}
&S(x,y,z) =  \int_\L \langle e_\l, \pi_\l(x,y,z)e_\l\rangle |\l| d\l \\
&=  - \int_\L  \bold 1_{[-1,0]}(\l) \langle{\bf 1}_{[-1,0]} , \pi_\l(x,y,z){\bf 1}_{[-1,0]} \rangle  \l d\l   + 
 \int_\L  \bold 1_{[0,1]}(\l) \langle  \bold 1_{\left[\frac{1}{\l} - 1, \frac{1}{\l}\right]} , \pi_\l(x,y,z) \bold 1_{\left[\frac{1}{\l} - 1, \frac{1}{\l}\right]} \rangle  \l d\l \\
 &= - \int_\L  \bold 1_{[-1,0]}(\l)  \int_\R \  e^{-2\pi i \l z}  e^{2\pi i \l yt} {\bf 1}_{[-1,0]}  {\bf 1}_{[-1,0]}(t-x)dt  \l d\l  \\
&\hspace{1in}+ \int_\L  \bold 1_{[0,1]}(\l)  \int_\R e^{-2\pi i \l z}  e^{2\pi i \l yt}   {\bold 1}_{\left[\frac{1}{\l} - 1, \frac{1}{\l}\right]}{\bold 1}_{\left[\frac{1}{\l} - 1, \frac{1}{\l}\right]}(t-x)dt  \l d\l  \\ 
&= - \int_\L   \bold 1_{[-1,0]}(\l)  \  e^{-2\pi i \l z} \left( \int_\R \ e^{2\pi i \l yt} {\bf 1}_{J_x}(t) dt\right)  \l d\l  \\
&\hspace{1in}+ \int_\L  \bold 1_{[0,1]}(\l)  e^{-2\pi i \l z}\left( \int_\R  e^{2\pi i \l yt}   {\bold 1}_{I_{x,\l}}(t) dt\right)  \l d\l  \\ 
&= \int_\L \ F_{x,y}(\l) e^{-2\pi i \l z} d\l +  \int_\L \ G_{x,y}(\l) e^{-2\pi i \l z} d\l.
\end{aligned}
$$
The explicit expression for $G$ is now an elementary calculation.
  
\end{proof}

We conclude this section with  a necessary and sufficient condition  for the generator of a Heisenberg orthonormal basis for any arbitrary shift-invariant spaces. Let  $g \in L^2(\L\times\R)$ and define the closed subspace $\mathcal S(g,\a,\b)$ of $L^2(\L\times\R)$ by
$$
\mathcal S(g,\a,\b) = \overline{sp} \bigl(\widehat{\mathcal T}(g, \a, \b)\bigr).
$$
For each $(\l,t) \in \L\times\R$ put
\begin{align}
\Theta^g_k(\l, t):= \sum_{l'\in\frac{1}{\b}\Bbb Z, l''\in \Bbb Z} g\left({\l - l''} , \frac{t-l'}{\l-l''}-k\right) \overline{g}\left({\l - l''} , \frac{t-l'}{\l-l''}\right) .
\end{align}

Then we have the following 

\begin{thm} \label{ONcharacter}   $\widehat{\mathcal T}(g, \a, \b)$ is an orthonormal basis for $\mathcal S(g,\a,\b)$ if and only if 
$$ \Theta^g_k (\l, t)= \delta_k\quad a.e.~ (\l, t).$$
\end{thm}

 

\begin{proof} For convenience we consider the case  $\a=\b=1$; the proof for general $\a$ and $\b$  can be adapted. For each $\gamma  = (k,l,m) \in \Gamma_{1,1}$ the function 
$$(\l,t) \mapsto  e^{2\pi i\l m} e^{-2\pi i \l l t} g(\l, t-k) \overline{g}(\l,t)| \l |$$  is absolutely integrable and we can apply periodization and Fubini's theorem  to calculate
$$
\begin{aligned}
\langle \hat T_\gamma g, g\rangle &= \int_\L  \int_\R  e^{2\pi i\l m} e^{-2\pi i \l l t} g(\l, t-k) \overline{g}(\l,t)| \l |   dt d\l \\
&= \int_\L  \int_\R e^{2\pi i\l m} e^{-2\pi i   l t} g(\l, t/\l-k) \overline{g}(\l , t/\l)dt d\l\\
&= \int_\L  e^{2\pi i\l m} \ \sum_{l'\in \Z} \  \int_0^1 e^{-2\pi i   l t} g(\l,(t-l')/\l-k) \overline{g}(\l,(t-l')/\l)~ dt d\l\\
&= \int_0^1 \int_0^1  e^{2\pi i\l m} e^{-2\pi i   l t}  \ \sum_{l''} \sum_{l'\in \Z} g\left(\l - l'', \frac{t-l'}{\l-l''}-k\right) \overline{g}\left(\l - l'', \frac{t-l'}{\l-l''}\right) dt d\l\\
&= \int_0^1 \int_0^1\ e^{2\pi i\l m} e^{-2\pi i   l t} \   \Theta^g_k (\l, t) dt d\l
\end{aligned}
$$
Suppose that $\widehat{\mathcal T}(g, \a, \b)$ is  an orthonormal basis for $\mathcal S(g,\a,\b)$. Note that $\Theta^g_k$ is a $(1,1)$-periodic integrable function on $\T\times \T$. If $k \ne 0$, then $\widehat{\Theta^g_k}(m,l) = 0 $ for all integers $m$ and $l$, and hence $\Theta^g_k \equiv 0$. If $k  = 0$, then $\widehat{\Theta^g_0}(m,l) = 0 $ holds for all $(m,l) \ne 0$ while $\widehat{\Theta^g_0}(0,0) = 1 $. Hence $\Theta^g_0\equiv 1 $.

On the other hand, if $ \Theta^g_k (\l, t)= \delta_k\quad a.e.~ (\l, t)$, then the above reasoning can be reversed to show that the system $\widehat{\mathcal T}(g, \a, \b)$  is orthonormal.  

\end{proof}

\vspace{.5cm}

 Bradley Currey,
 {Department of Mathematics and Computer
Science, Saint Louis University, St. Louis, MO 63103}\\
 { \footnotesize{E-mail address: \texttt{{ curreybn@slu.edu}}}\\}

  Azita Mayeli,
  {Mathematics Department, City College of Technology, City University of New York,  New York, USA }\\
   \footnotesize{E-mail address: \texttt{{amayeli@citytech.cuny.edu}}}\\


\begin{thebibliography}{99}
\bibitem{CM08}
{B. ~Currey, A.~ Mayeli,} {\it Gabor Fields and Wavelet Sets for the Heisenberg Group}, preprint. 
\bibitem{F}
{ H. ~F\"uhr,} {\it Abstract Harmonic Analysis of Continuous Wavelet
Transforms,} {  Lect. Notes in Math.} {\bf 1863} \ (2005), Springer

\bibitem{FG} {\sc H. ~F\"uhr, K. Gr\"ochenig,}
{\it Sampling theorems on locally compact groups from oscillation estimates}, {
Math. Z.} {\bf 255}  \ (2007),  177--194.
\bibitem{M08} {\sc A. ~Mayeli,}
{\it Shannon multiresolution analysis on the Heisenberg group}, {J.
Math. Anal. Appl.} {\bf 348}  \ (2008),  No. 2, 671-684.
\bibitem{Pesen98} I. ~ Pesenson, {\it Sampling of Paley-Wiener functions on stratified groups}, J. Fourier Anal. Appl. 4 (1998) 271 -- 281.
\end{thebibliography}
\end{document}